\documentstyle[11pt]{article}
\font\fr=eufm10   scaled \magstep 1   
\font\Bbb=msbm10 scaled \magstep 1   
\font\bfpeq=cmbx7 scaled \magstep 1
\def\zeta{Z\! \! \! Z}                      

\newcommand{\fin}{\hfill$\Box$}
\newcommand{\com}{\makebox[7pt]{\raisebox{1.5pt}{\tiny
$\circ$}}}

\begin{document}
\date{}
\section*{The number of functionally independent invariants of a
pseudo--Riemannian metric}

\vspace{4mm}

\hfill
\parbox{11.2cm}{
    {J Mu\~noz Masqu\'e\dag\, and Antonio Vald\'es\ddag} \\
    {\footnotesize \dag CSIC--IEC, C/Serrano 144, 28006--Madrid, Spain}\\
    {\footnotesize \ddag UNED, Departamento de Matem\'aticas Fundamentales,
    C/Senda del Rey s/n, 28040--Madrid, Spain}\\[4mm]

    {\bfpeq Short title:} {\footnotesize Invariants of a pseudo--Riemannian
    metric}\\[4mm]

    {\footnotesize PACS numbers: $02.40.+m \;\;\; 04.20.Cv$}\\[4mm]

    {\bfpeq Abstract.} {\footnotesize The number of functionally independent
    scalar invariants of arbitrary order of a generic pseudo--Riemannian
    metric on an $n$--dimensional manifold is determined.}}

\vspace{6mm}

\noindent {\bf 1. Introduction}

\vspace{5mm}

\noindent The goal of this work is to determine the number $i_{n,r}$ of
functionally independent differential invariants of order $r$ of a
{\em generic} pseudo--Riemannian metric $g$ on an $n$--dimensional
manifold $N$. The results are as follows: for every $n\geq 1$,
$i_{n,0} = i_{n,1} = 0$; for every $r\geq 0$, $i_{1,r} = 0$; $i_{2,2}
= 1$, and for every $r\geq 3$, $i_{2,r} = \frac{1}{2}(r+1)(r-2)$; finally,
$$i_{n,r} = n + \frac{(r-1)n^2-(r+1)n}{2(r+1)}{n+r\choose r}\,\,,\,\,\,
\mbox{\rm for every}\;\;\, n\geq 3\,,\,r\geq 2\,.$$
The theory of metric invariants is classic in both General Relativity and
Riemannian Geometry. The standard approach to this topic relies on the
definition of an invariant as a polynomial in the $g_{ij}$'s, their
partial derivatives up to a certain order, say $\partial^{|\alpha|}g_{ij}/
\partial x^{\alpha}$, $|\alpha|\leq r$, and in $[\hbox{\rm det}(g_{ij})]
^{-1}$, which is ``natural'' under diffeomorphisms (for example, see [1]).
For scalar invariants the above definition does not
allow one to pose the question of how many functionally
independent invariants there are for each order since some of the funtional
relationships among invariants may be out of the ring that the above
definition prescribes and furthermore, since that ring is not complete
standard tools of analysis (as involutiveness, Frobenius theorem, etc.)
cannot be applied either. From this point of view, the enumeration of
the scalars constructed from the Riemann tensor of the Levi--Civit\`a
connection of a pseudo--Riemannian metric by means of covariant
differentiation, tensor products and contractions has been discussed
in some recent papers: in [5] the number of
independent homogeneous scalar monomials of each order
and degree up to order $12$ in the derivatives of the metric is determined
and in [10] the same number is determined up to order $14$.
Apart from the interest and complexity of these results specially in
relation to the so called Weyl invariants (cf. [4]) for the field
theory it is clear that the determination of the number $i_{n,r}$ is
the most relevant fact since it provides the number of essentially
different Diff$(N)$--invariant Langrians of arbitrary order that exist
in G.R. It seems thus natural to found the theory
on the jet bundle notion of an invariant (cf. [7]) which avoids the
aforementioned difficulties of the polynomial notion
and translates the naturality condition into an
authentic condition of invariance under the action of the group of
diffeomorphisms of $N$ on an appropriate jet bundle.

The plan of this paper is the following: In section $2$ we introduce
the notion of a metric invariant as well as that of an invariant
Lagrangian density although for an oriented ground manifold $N$ the
latter is reduced to the former since the bundle of metrics
over $N$ is endowed with a canonical invariant zero order Lagrangian
density, so that the emphasis is put on scalar invariants. The notion of
invariance is related to a specific representation of the vector fields
of $N$ into vector fields of the $r$--jet bundle of metrics. Section $3$
contains the explicit determination of this representation and its formulas
are used throughout the paper. In section $4$  we prove that on a dense
open subset of the $r$--jet bundle the metric invariants coincide with the
ring of first integrals of an involutive distribution which is obtained
by linearizing the basic representation by means of a homomorphism of
vector bundles, $\Phi^r$. The number of invariants $i_{n,r}$, is thus
equivalent to know the rank of $\Phi^r$. Sections $5, 6$ and $7$ are
devoted to this aim distinguishing the different cases that appear according
to the values of the order $r$ of the jet bundle that we are considering
and the dimension $n$ of the ground manifold. Finally, section $8$
contains the calculation of $i_{n,r}$ and the comparison of $i_{n,2}$ with
the standard procedure (cf. [9]) in order to generate the second order
metric invariants.

\vspace{8mm}

\noindent {\bf 2. The notion of a metric invariant}

\vspace{5mm}

\noindent Let $N$ be an $n$--dimensional differentiable manifold.
Given an integer $0\leq n^+\leq n$, we shall denote by $p:{\cal M} =
{\cal M}_{n^+}(N)\rightarrow N$ the bundle of pseudo--Riemannian metrics
on $N$ of signature $(n^+,n^-)$, $n^- = n - n^+$ ({\em i.e.},
the global sections of $p$ are exactly the pseudo--Riemannian
metrics on $N$ of signature $(n^+,n^-)$ at each point). Let $p_r:J^r({\cal M})
\rightarrow N$ be the $r$--jet bundle of local sections of $p$. The
$r$--jet at a point $x\in N$ of a metric $g$ of ${\cal M}$ will be denoted
by $j^r_x(g)$. For every $r\geq s$, we also have a natural projection
$p_{rs}:J^r({\cal M})\rightarrow J^s({\cal M})$, $p_{rs}(j^r_x g) = j^s_x g$.
Let $(U; x_1,\ldots,x_n)$ be an open coordinate domain of $N$ and let
$\alpha=(\alpha_1,\ldots,\alpha_n)$ be a multi--index of non--negative
integers. We set $|\alpha|=\sum_i\alpha_i.$ The family of functions
$(x_i\com p_r,y^{jk}_{\alpha}),\;j\leq k, |\alpha|\leq r$,
defined by $y^{jk}_\alpha(j^r_xg)=
\left(\partial^{|\alpha|}g_{jk}/\partial x^{\alpha}\right)(x)$,
where $g_{jk}=g(\partial/\partial x_j,
\partial/\partial x_k)$, constitutes a coordinate chart on $p_r^{-1}U =
J^r(p^{-1}U)$. We shall simply write $y_{jk}$ instead of $y^{jk}_0$. Note
that the functions $(x_i\com p, y_{jk})$, $1\leq i\leq n$, $1\leq j\leq
k\leq n$, are a coordinate system on $p^{-1}(U)$. We shall also set
$y^{jk}_{\alpha} = y^{kj}_{\alpha}$, for $j>k$.

Let $f\colon N\rightarrow N'$ be a diffeomorphism.
We shall denote by $\overline{f}\colon {\cal M}\rightarrow {\cal M}'$,
${\cal M}'={\cal M}_{n^+}(N')$,
the natural lift of $f$ to the bundles of pseudo--Riemannian metrics;
{\em i.e.}, $\overline{f}(g_x) = (f^{-1})^* g_x$. Hence $p'\com
\overline{f} = f\com p$.

The diffeomorphism $\overline{f}\colon {\cal M}\rightarrow {\cal M}'$ has a
natural extension to jet bundles, $J^r(f)\colon J^r({\cal M})\rightarrow
J^r({\cal M}')$, defined as follows:
$J^r(f)\left(j^r_xg\right)= j^r_{f(x)}(\overline{f}\com g\com
f^{-1})$.

Given a vector field $X\in \mbox{\fr X}(N)$,  we shall denote by
$\overline{X}^r$ the natural lift of $X$ to $J^r({\cal M})$. For $r=0$ we
shall simply write $\overline{X}$ instead of $\overline{X}^0$. Note
that $\overline{X}$ is the natural lift to the bundle ${\cal M}$ of
pseudo--Riemannian metrics of the vector field $X$.
If $\varphi_t$ is the local flow of $X$, then $J^r(\varphi_t)$ is the
local flow of $\overline{X}^r$. Hence $\overline{X}^r$ is projectable
onto $\overline{X}^{r-1}$, and $\overline{X}$ is projectable onto $X$.
The mapping $X\mapsto \overline{X}^r$, is a \mbox{\Bbb R}--linear injection
and for every $X,Y\in \mbox{\fr X}(N)$,
\begin{equation}
\left[\overline{X}^r, \overline{Y}^r\right] = \overline{[X,Y]}^r.
\end{equation}
Hence we have a faithful representation of
\mbox{\fr X}$(N)$ into \mbox{\fr X}$\left(J^r({\cal M})\right)$.

\vspace{2mm}

\noindent {\em Definition 1.} A function $F\in C^\infty(J^r({\cal M}))$
(which may be only defined on an open subset) is said to be a
{\em metric differential invariant} of order $r$ if for every
$X\in \mbox{\fr X}(N)$, $\overline{X}^rF=0$.

\vspace{2mm}

\noindent {\em Definition 2.} A function $F\in C^\infty(J^r({\cal M}))$
is said to be a {\em metric invariant} of order $r$ if for every
diffeomorphism $f:N\rightarrow N$, $F\com J^r(f) = F$.

\vspace{2mm}

\noindent {\em Remark 3.} Metric invariants are a subring of the ring
of metric differential invariants. In fact, the set of vector fields
on $N$ with compact support, \mbox{\fr X}$_c(N)$, is a dense ideal of
\mbox{\fr X}$(N)$ with respect to the $C^\infty$ topology and hence
a function $F\in C^\infty(J^r({\cal M}))$ is a metric differential invariant
if and only if for every $X\in \mbox{\fr X}_c(N)$, $\overline{X}^r F = 0$.
This is equivalent to saying that for every $t\in \mbox{\Bbb R}$, one
has $F\com J^r(\phi_t) = F$, $\phi_t$ being the one--parameter group
generated by $X$, and the last equation evidently holding for a metric
invariant.

\vspace{2mm}

\noindent {\em Example 4.} Let $\nabla$ be the Levi--Civit\`a connection
of a pseudo--Riemannian metric $g$ of ${\cal M}$, and $R$ the curvature tensor.
Since $R$ is of type $(1,3)$, for every $r\in \mbox{\Bbb N}$,
$\nabla^{2r}R$ is a tensor field of type $(1,2r+3)$. Let us choose a sequence
of $r+1$ covariant indices $1\leq i_0<...<i_r\leq 2r+3$, and let us
apply to them the isomorphism $g^{\sharp}:T_x^*(N)\rightarrow T_x(N)$,
thus obtaining a tensor field ${}^{g^{\sharp}}\left(\nabla^{2r}R\right)
^{i_0,...,i_r}$ of type $(r+2,r+2)$. If further we choose a permutation
$j_1,...,j_{r+2}$ of its covariant indices then we can obtain an scalar
by simply setting $S_g = c^1_{j_1}\cdots c^{r+2}_{j_{r+2}}\left(
{}^{g^{\sharp}}\left(\nabla^{2r}R\right)^{i_0,...,i_r}\right)$, where
$c^i_j$ stands for the contraction of the $i$--th contravariant index with
the $j$--th covariant index.
The value of $S_g$ at a point $x\in N$ only depends on
$j^{2r+2}_x(g)$, since the local coefficients $\Gamma^i_{jk}(x)$ of
$\nabla$ only depend on $j^1_x(g)$,
and $R_x$ only depends on $j^2_x(g)$ (cf. [6], IV.2.4 and III.7.6).
Accordingly, we can define a function $F\in C^\infty(J^{2r+2}({\cal M}))$ by
imposing that $F(j^{2r+2}_xg) = S_g(x)$, and this function is an invariant.
In fact, if $f:N\rightarrow N$ is a diffeomorphism and we set $\bar{g} =
(f^{-1})^* g = \overline{f}\com g\com f^{-1}$, then the Levi--Civit\`a
connection of $\bar{g}$ is the linear connection $\overline{\nabla}$ given
by $\overline{\nabla}_X Y = f\cdot\left(\nabla_{f^{-1}\cdot X} f^{-1}\cdot
Y\right)$, as follows from Koszul's formula ([6], IV.2.3), and consequently
for every $r\in \mbox{\Bbb N}$, the tensor fields $\nabla^{2r}R$
and $\overline{\nabla}^{2r}\bar{R}$ are $f$--related (cf.
[6], VI.1.2); {\em i.e.}, for every system of vector fields $X_1,...,X_{2r+1}
\in \mbox{\fr X}(N)$, and every point $x\in N$, $f_*\left(\nabla^{2r}R\right)
\left((X_1)_x,...,(X_{2r+3})_x\right) = \left(\overline{\nabla}^{2r}\bar{R}
\right)\left((f\cdot X_1)_{f(x)},...,(f\cdot X_{2r+3})_{f(x)}\right)$, or
else  $\left(\nabla^{2r}R\right)\left(X_1,...,X_{2r+3}\right) =
f^{-1}\cdot\left(\overline{\nabla}^{2r}\bar{R}\right)\left(f\cdot X_1,
...,f\cdot X_{2r+3}\right)$. Hence $S_g(x) = S_{\bar{g}}(f(x))$, and this
means $F(j^{2r+2}_x g)$ $= F\left(J^{2r+2}(f)(j^{2r+2}_x g)\right)$.

\vspace{2mm}

\noindent {\em Definition 5.} An $r^{th}$ order {\em Lagrangian density}
is a horizontal $n$--form $\Omega_n$ on $J^r({\cal M})$. An $r^{th}$ order
Lagrangian density is said to be {\em invariant} if for every $X\in
\hbox{\fr X}(N)$, $L_{\overline{X}^r}\Omega_n=0$.

\vspace{2mm}

\noindent {\em Remark 6.} As $\Omega_n$ is horizontal, locally there
exists a function ${\cal L}\in C^\infty(J^r{\cal M})$, such that $\Omega_n =
{\cal L}\,\,dx_1\wedge...\wedge dx_n$. Below we shall see that by
introducing the factor $\sqrt{(-1)^{n^-}det(y_{ij})}$ in $dx_1\wedge...
\wedge dx_n$ we obtain a globally defined invariant Lagrangian density,
thus reducing the problem of determining the invariant Lagrangian
densities to that of the scalar invariants. Note that in this case
${\cal L}$ should be substituted by $F= {\cal L}/
\sqrt{(-1)^{n^-}det(y_{ij})}$.

\vspace{2mm}

\noindent {\em Proposition 7.} Assume $N$ is oriented. Then the bundle of
metrics ${\cal M}$ is endowed with a canonical invariant zero order
Lagrangian density $\omega_n$, uniquely defined by the following
condition: if $X_1,...,X_n$ is an orthonormal basis for a metric
$g$ of ${\cal M}$, defined on an open subset $U\subset {\cal M}$,
which belongs to the orientation of $N$, then for every $x\in U$,
$(\omega_n)_{g_x}(\overline{X}_1,...,\overline{X}_n) = 1$.
Accordingly, every Lagrangian density $\Omega_n$ on $J^r({\cal M})$
can be uniquely written as $\Omega_n = F\,\,\omega_n$, $F\in
C^\infty\left(J^r({\cal M})\right)$, and $\Omega_n$ is invariant if and only if
$F$ is a metric differential invariant.

\vspace{2mm}

\noindent {\em Proof.} Since $\omega_n$ must be a horizontal $n$--form
it is clear that the condition in the statement uniquely determines the
desired form. Moreover, we can define a horizontal $n$--form on ${\cal M}$
by setting for every $Y_1,...,Y_n\in T_{g_x}({\cal M})$, $(\omega_n)_{g_x}
(Y_1,...,Y_n) = v_{g_x}(p_*Y_1,...,p_*Y_n)$, where $v_{g_x}$ is the
Riemannian volume associated with $g_x$. Since $\overline{X}_i$ is
$p$--projectable onto $X_i$, we have $(\omega_n)_{g_x}(\overline{X}_1,
...,\overline{X}_n) = v_{g_x}(X_1,...,X_n)$, thus proving that $\omega_n$
satisfies the above condition.

A basis $X_1,...,X_n$ for $T_x(N)$ is said to be orthonormal for the
metric $g_x$ if: $g(X_i,X_j) = \delta_{ij}$ for either $1\leq i\leq n^+$
or $1\leq j\leq n^+$; $g(X_i,X_j) = - \delta_{ij}$ for
$n^+ + 1\leq i,j\leq n^+ + n^-$; in other words, the matrix of $g_x$ must be
$$\left(
\begin{array}{cc}
I_{n^+} & 0 \\
0 & -I_{n^-}
\end{array}
\right).
$$
Hence locally we have $\omega_n = \sqrt{(-1)^{n^-}\hbox{\rm det}(y_{ij})}
dx_1\wedge...\wedge dx_n$. Also note that $\omega_n$ cannot be considered
as the volume element associated to the canonical metric $G =
\sum_{i\leq j}y_{ij} dx_i\otimes dx_j$ on ${\cal M}$, since $G$ is singular!

We shall now prove that $\omega_n$ is invariant. Given a diffeomorphism
$f$ of $N$, with the above notations we have: $\left(\overline{f}^*\omega_n
\right)_{g_x}(Y_1,...,Y_n)=$ $(\omega_n)_{\overline{f}(g_x)}(\overline{f}_*
Y_1,...,\overline{f}_*Y_n)=$ $v_{\overline{f}(g_x)}(p_*\overline{f}_*Y_1,
...,p_*\overline{f}_*Y_n)$, and since $p_*\com \overline{f}_*=f_*\com p_*$,
we obtain\newline
\hspace*{2.4cm} $\left(\overline{f}^*\omega_n\right)_{g_x}(Y_1,...,Y_n)=$
$\left(f^* v_{\overline{f}(g_x)}\right)(p_*Y_1,...,p_*Y_n)=$\newline
$\left(f^* v_{(f^{-1})^*(g_x)}\right)(p_*Y_1,...,p_*Y_n)=$
$f^*(f^{-1})^* v_{g_x}(p_*Y_1,...,p_*Y_n)=$ $\left(\omega_n\right)_{g_x}
(Y_1,...,Y_n)$.

\vspace{8mm}

\noindent {\bf 3. Local expression of the basic representation}

\vspace{5mm}

\noindent {\em Proposition 8.} Let $X = \sum_i u_i (\partial/\partial x_i)$,
$u_i\in C^\infty(U)$, $1\leq i\leq n$, be the local expression of a vector
field $X\in \mbox{\fr X}(N)$ on an open coordinate domain $(U;x_1,...,x_n)$
of $N$. The local expression of the lifting of $X$ to the bundle of
pseudo--Riemannian metrics, $\overline{X}\in \mbox{\fr X}({\cal M})$, in the
induced coordinate system $(p^{-1}(U); x_i, y_{jk})$, $1\leq i\leq n$,
$1\leq j\leq k\leq n$, is given by
\begin{equation}
\overline{X} = \sum_i u_i \frac{\partial}{\partial x_i} + \sum_{i\leq
j}v_{ij}\frac{\partial}{\partial y_{ij}} \hspace*{1cm} v_{ij} = -\sum_h
\frac{\partial u_h}{\partial x_i}y_{hj} - \sum_h \frac{\partial
u_h}{\partial x_j}y_{ih}\,.
\end{equation}

\vspace{2mm}

\noindent {\em Proof.} First, note that $v_{ij}$ is symmetric with respect
to the indices $i,j$, so we shall also write $v_{ij} = v_{ji}$, for $i>j$.
As is well--known, the lift $\overline{X}$ is the
unique vector field on ${\cal M}$ which is $p$--projectable onto $X$ and
leaves the ``canonical metric'', $G = \sum_{i\leq j}y_{ij} dx_i
\otimes dx_j$, on the manifold ${\cal M}$, invariant; {\em i.e.},
$\overline{X} = \sum_i u_i(\partial/\partial x_i) + \sum_{i\leq j}v_{ij}
(\partial/\partial y_{ij})$, for some functions $v_{ij}\in C^\infty(p^{-1}
U)$, and $L_{\overline{X}}G = 0$. Hence
$$L_{\overline{X}^r}G= \sum_{i\leq j}\left[v_{ij}dx_i\otimes dx_j+ y_{ij}
du_i\otimes dx_j + y_{ij}dx_i\otimes du_j\right] = 0\,,$$
and this equation completely determines the unknown functions.
\fin

\vspace{2mm}

{}From the general formulas for the prolongation of vector fields
by infinitesimal contact transformations ({\em e.g.}, see [8]) we then
obtain the local expression for $\overline{X}^r$; more precisely,
\begin{equation}
\overline{X}^r=\sum_iu_i\frac{\partial}{\partial x_i}-\sum_h\sum_{i\leq
j}\sum_{|\alpha|=0}^r\left\{\sum_{\beta\leq \alpha}{\alpha\choose\beta}
\left[\frac{\partial^{|\beta|+1}u_h}{\partial x^{\beta+(i)}}
y^{hj}_{\alpha-\beta}+\frac{\partial^{|\beta|+1}u_h}{\partial x^{\beta+(j)}}
y^{ih}_{\alpha-\beta}\right]\;+\right.
\end{equation}
$$\hspace*{3.5cm}\left. \sum_{0<\beta\leq \alpha}{\alpha\choose\beta}
\frac{\partial^{|\beta|}u_h}{\partial x^{\beta}}y^{ij}_{\alpha-\beta+
(h)}\right\}\frac{\partial}{\partial y^{ij}_{\alpha}}\,,$$
where $(i)$ stands for the multi--index $(i) = (0,...,\stackrel{i}{1},...,0)$.
The above equations can be obtained by either imposing that: $1^{st})\;
\overline{X}^r$ is $p_r$--projectable onto $X$, and $2^{nd})
\; \overline{X}^r$ leaves
the generalized contact differential system ${\cal C}$
spanned by the one--forms on $J^r({\cal M})$, $\theta^{ij}_{\alpha} =
dy^{ij}_{\alpha}-\sum_k y^{ij}_{\alpha+(k)} dx_k$, $i,j=1,...,n$,
$|\alpha|<r$ invariant; {\em i.e.}, $L_{\overline{X}^r}{\cal C}\subset
{\cal C}$, or by simply calculating the infinitesimal generator associated
to $J^r(\phi_t)$, $\phi_t$ being the local flow of $X$.

\vspace{2mm}

\noindent {\em Example 9.} For $r=1$, the above formula reads as follows:
\begin{equation}
\overline{X}^1= \sum_i u_i \frac{\partial}{\partial x_i}
-\sum_h\sum_{i\leq j}\left(\frac{\partial u_h}{\partial x_i}y_{hj}+
\frac{\partial u_h}{\partial x_j}y_{ih}\right)\frac{\partial}{\partial
y_{ij}}\; -
\end{equation}
$$- \sum_{h,k}\sum_{i\leq j}\left(\frac{\partial^2 u_h}{\partial x_i\partial
x_k}y_{hj} + \frac{\partial^2 u_h}{\partial x_j\partial
x_k}y_{ih} + \frac{\partial u_h}{\partial x_i}y^{hj}_k +
\frac{\partial u_h}{\partial x_j}y^{ih}_k +
\frac{\partial u_h}{\partial x_k}y^{ij}_h\right)\frac{\partial}{\partial
y^{ij}_k}.$$

\vspace{8mm}

\noindent {\bf 4. The fundamental distribution}

\vspace{5mm}

\noindent {\em Theorem 10.} With the above hypotheses and notations we have:
\begin{description}
\item[(i)] $\overline{X}^r_{j^r_x g}$ only depends on $j^{r+1}_x(X)$.
\item[(ii)] There exists a unique homomorphism of vector bundles over
$J^r({\cal M})$,
$$\Phi^r:p_r^* J^{r+1}(TN) \longrightarrow T(J^r{\cal M})\,,$$
such that for every $X\in \mbox{\fr X}(N)$, $j^r_x g\in J^r({\cal M})$,
$$\Phi^r(j^r_x g , j^{r+1}_x X) = \overline{X}^r_{j^r_x g}\,.$$
\item[(iii)] On a dense open subset ${\cal O}^r\subset J^r({\cal M})$,
the image of $\Phi^r$ defines an involutive distribution $\hbox{\fr D}^r$,
such that for every $j^r_x g\in {\cal O}^r$, $X\in \hbox{\fr X}(N)$,
$$\hbox{\fr D}^r_{j^r_x g} = \left\{\overline{X}^r_{j^r_x g}\;,\;X\in
\hbox{\fr X}(N)\right\} \subset T_{j^r_x g}(J^r{\cal M})\,.$$
\item[(iv)] A function $F\in C^\infty({\cal O}^r)$ is a metric differential
invariant if and only if $F$ is a first integral of $\hbox{\fr D}^r$.
\end{description}

\vspace{2mm}

\noindent {\em Proof.} Evaluating $\overline{X}^r$ at the point
$j^r_x g$, from formula (3) we obtain:
$$\overline{X}^r_{j^r_x g}(y^{ij}_{\alpha}) =$$
$$-\sum_h\left\{\sum_{\beta\leq \alpha}{\alpha\choose\beta}
\left[\frac{\partial^{|\beta|+1}u_h}{\partial x^{\beta+(i)}}(x)
\frac{\partial^{|\alpha-\beta|} g_{hj}}{\partial x^{\alpha-\beta}}(x)
+\frac{\partial^{|\beta|+1}u_h}{\partial x^{\beta+(j)}}(x)
\frac{\partial^{|\alpha-\beta|}g_{ih}}{\partial x^{\alpha-\beta}}(x)\right]\;
+\right.$$
$$\hspace*{7cm}\left. \sum_{0<\beta\leq \alpha}{\alpha\choose\beta}
\frac{\partial^{|\beta|}u_h}{\partial x^{\beta}}(x)\frac{\partial^
{|\alpha-\beta|+1}g_{ij}}{\partial x^{\alpha-\beta+(h)}}(x)\right\}\,,
$$
thus proving (i). Accordingly, we can define a unique map $\Phi^r$,
by setting: $\Phi^r(j^r_x g , j^{r+1}_x X) = \overline{X}^r_{j^r_x g}$.
Since the map $X\mapsto \overline{X}^r$ is \hbox{\Bbb R}--linear, it is
clear that $\Phi^r$ is a homomorphism of vector bundles as stated in (ii).

Let us define the subset ${\cal O}^r$ as follows:
a point $j^r_x g$ belongs to ${\cal O}^r$ if and only if it has
a neighbourhood ${\cal N}_{j^r_x g}$ such that the rank of
$\Phi^r_{|{\cal N}_{j^r_x g}}$ is constant. From the very definition, ${\cal
O}^r$ is an open subset and the rank of $\Phi^r_{|{\cal O}^r}$ is locally
constant. Next we prove that ${\cal O}^r$ is dense in $J^r({\cal M})$.
Let ${\cal U}$ be a non--empty open subset of $J^r({\cal M})$. Since
the rank of
$\hbox{\fr D}^r$ only takes a finite number of values, there exists a point
$j^r_x g\in {\cal U}$, such that for every $j^r_{x'} g'\in {\cal U}$,
rk.$\hbox{\fr D}^r_{j^r_{x'} g'}\leq$ rk.$\hbox{\fr D}^r_{j^r_{x} g}$, and
since the rank of a homomorphism of vector bundles is a lower
semicontinuous function, $j^r_x g\in {\cal O}^r$. In order to prove that
$\hbox{\fr D}^r$ is involutive we proceed as follows. Given a point
$j^r_x g\in {\cal O}^r$, let $X_1,...,X_k$ be vector fields on $N$ such
that $\left(\overline{X}_1\right)^r_{j^r_x g},...,
\left(\overline{X}_k\right)^r_{j^r_x g}$ is a basis for $\hbox{\fr D}^r_{
j^r_x g}$. Then, there exists an open neighbourhood ${\cal N}_{j^r_x g}$
such that $\overline{X}^r_1,...,\overline{X}^r_k$ is a basis of $\hbox{
\fr D}^r_{j^r_{x'} g'}$, for every $j^r_{x'} g'\in {\cal N}_{j^r_x g}$.
Accordingly any two vector fields $\xi$, $\xi'$, belonging to
$\hbox{\fr D}^r_{|{\cal N}_{j^r_x g}}$ can be written as $\xi = \sum_i f_i
\left(\overline{X}^r_i\right)$, $\xi' = \sum_j f'_j \left(\overline{X}^r
_j\right)$, and from formula (1) we obtain
$$[\xi, \xi'] = \sum_{i,j=1}^k \left\{f_i\,\overline{X}^r_i(f'_j)\,
\overline{X}^r_j - f'_j\,\overline{X}^r_j(f_i)\,\overline{X}^r_i
+ f_if'_j\overline{[X_i,X_j]}^r\right\}\,,$$
thus showing that $[\xi, \xi']$ also belongs to $\hbox{\fr D}^r$, and
finishing the proof of (iii). Part (iv) follows directly from the
definitions.
\fin

\vspace{4mm}

\noindent {\em Corollary 11.} On a neighbourhood of each point $j^r_x g\in
{\cal O}^r$, the number of functionally independent metric differential
invariants is\,\, dim$J^r({\cal M}) - $ rk.$\Phi^r_{j^r_x g}\,$.

\vspace{2mm}

\noindent {\em Proof.} This follows from Theorem 7 and the Frobenius theorem.
\fin

\vspace{2mm}

Our next goal is to determine the rank of $\Phi^r$. In doing this we shall
use normal coordinates which will always be assumed to be metric ({\em i.e.},
associated with an orthonormal frame) and defined on a convex open
neighbourhood of a given point $x\in N$ ([6], III.\S8, IV.\S3). The expansion
of the metric in a normal coordinate system starts as follows (cf. [3]):
\begin{equation}
g_{ij} = g_{ij}(x) + \frac{1}{6} \sum_{k,l = 1}^n
\left(R_{ilkj}(x)+ R_{jlki}(x)\right) x_k x_l\, + \,\, \cdots \,\,,
\end{equation}
where $R_{ijkl}$ are the components of the curvature tensor, {\em i.e.},
$$R_{ijkl} =
g\left(R(\partial/\partial x_k, \partial/\partial x_l)(\partial/\partial x_j)
\,,\,\partial/\partial x_i\right)\,.$$
Taking derivatives in (5), we obtain
\begin{equation}
\frac{\partial g_{ij}}{\partial x_k}(x) = 0 \hspace*{1cm} 1\leq i,j,k\leq n\,,
\end{equation}
and again taking derivatives,
\begin{equation}
\frac{\partial^2g_{ij}}{\partial x_k \partial x_l}(x) = \frac{1}{3}
\left(R_{iklj}(x)+R_{ilkj}(x)\right)\,.
\end{equation}

\vspace{8mm}

\noindent {\bf 5. The rank of $\Phi^1$}

\vspace{5mm}

\noindent {\em Theorem 12.} With the same notations as in Proposition 8,
$\Phi^1(j^{2}_x X) = \overline{X}^1_{j^1_x g} = 0$ if and only if
in a normal coordinate system $x_1,...,x_n$ centred at $x\in N$, the
following conditions hold true: For every $i,j,k=1,...,n$,
\begin{equation}
u_i(x)=0\hspace*{1cm}g_{ii}(x)\frac{\partial u_i}{\partial x_j}(x)+
g_{jj}(x)\frac{\partial u_j}{\partial x_i}(x)=0\hspace*{1cm}
\frac{\partial^2 u_i}{\partial x_j\,\partial x_k}(x)=0\,.
\end{equation}
Hence $\Phi^1:J^2_x(TN)\rightarrow T_{j^1_x g}(J^1{\cal M})$ is surjective at
each point $j^1_x g\in J^1({\cal M})$.

\vspace{2mm}

\noindent {\em Proof.} From formula (4) we obtain
\begin{equation}
u_i(x) = 0\hspace*{1cm} g_{ii}(x)\frac{\partial u_i}{\partial x_j}(x) +
g_{jj}(x)\frac{\partial u_j}{\partial x_i}(x) = 0\,,
\end{equation}
\begin{equation}
g_{jj}(x)\frac{\partial^2u_j}{\partial x_i \partial x_k}(x) +
g_{ii}(x)\frac{\partial^2u_i}{\partial x_j \partial x_k}(x) +
\sum_h\left(\frac{\partial u_h}{\partial x_i}(x)\frac{\partial g_{hj}}{
\partial x_k}(x)+\right.
\end{equation}
$$\hspace*{5.5cm} \left.\frac{\partial u_h}{\partial x_j}(x)
\frac{\partial g_{ih}}{\partial x_k}(x)+
\frac{\partial u_h}{\partial x_k}(x)
\frac{\partial g_{ij}}{\partial x_h}(x)\right) = 0\,.$$
By applying (6), equation (10) becomes
\begin{equation}
g_{jj}(x)\frac{\partial^2u_j}{\partial x_i \partial x_k}(x) +
g_{ii}(x)\frac{\partial^2u_i}{\partial x_j \partial x_k}(x) = 0\,.
\end{equation}
By permuting $(i\,k)\mapsto (k\,i)$ in (11), we have
\begin{equation}
g_{jj}(x)\frac{\partial^2u_j}{\partial x_k \partial x_i}(x) +
g_{kk}(x)\frac{\partial^2u_k}{\partial x_j \partial x_i}(x) = 0\,.
\end{equation}
Comparing (11) and (12), we obtain
\begin{equation}
g_{kk}(x)\frac{\partial^2u_k}{\partial x_i \partial x_j}(x) =
g_{ii}(x)\frac{\partial^2u_i}{\partial x_j \partial x_k}(x) \,.
\end{equation}
{}From (13) and (11), and again applying (13) after making the permutation
$(i\,j\,k)\mapsto (j\,i\,k)$, we obtain
$$\underline{g_{kk}(x)\frac{\partial^2u_k}{\partial x_i \partial x_j}(x)} =
g_{ii}(x)\frac{\partial^2u_i}{\partial x_j \partial x_k}(x) =
-g_{jj}(x)\frac{\partial^2u_j}{\partial x_i \partial x_k}(x)=
-\underline{g_{kk}(x)\frac{\partial^2u_k}{\partial x_j \partial x_i}(x)}.$$
Accordingly,
$$\frac{\partial^2u_k}{\partial x_i \partial x_j}(x) = 0\,,$$
thus finishing the proof of (8). Hence
$$\hbox{\rm rk.}\Phi^1 = \hbox{\rm dim}J^2_x(TN) - \frac{n(n-1)}{2} =
\frac{n(n^2+2n+3)}{2} = \hbox{\rm dim}T_{j^1_x g}(J^1{\cal M})\,.
\hspace*{1.2cm} \hbox{\fin}$$

\vspace{8mm}

\noindent {\bf 6. The rank of $\Phi^2$}

\vspace{5mm}

\noindent {\em Lemma 13.} With the same notations as in Proposition 8
and Theorem 12, $\Phi^2(j^{3}_x X) = \overline{X}^2_{j^2_x g} = 0$
if and only if in addition to equations (8) the following conditions hold true:
For every $i,j,k,l=1,...,n$,
\begin{equation}
\frac{\partial^3 u_i}{\partial
x_j \partial x_k \partial x_l}(x)= 0\,,
\end{equation}

\begin{equation}
\sum_{h=1}^n\left(\frac{\partial
u_h}{\partial x_i}(x) R_{hkjl}(x) + \frac{\partial u_h}{\partial x_j}(x)
R_{hlik}(x) + \frac{\partial u_h}{\partial x_k}(x) R_{hilj}(x) +\right.
\end{equation}
$$\hspace*{8cm}\left.\frac{\partial u_h}{\partial x_l}(x) R_{hjki}(x)
\right)=0.$$

\vspace{2mm}

\noindent {\em Proof.} It follows from formula (3) that
$\overline{X}^2_{j^2_x g} = 0$ if and only if equations (8) hold and
furthermore for every $i,j,k,l=1,...,n$,
\begin{equation}
g_{ii}(x)\frac{\partial^3 u_i}{\partial x_j \partial x_k \partial x_l}(x)+
g_{jj}(x)\frac{\partial^3 u_j}{\partial x_i \partial x_k \partial x_l}(x)+
\lambda_{ijkl}=0\,,
\end{equation}
where we have set:
$$\lambda_{ijkl}=
\sum_{h=1}^n\left(\frac{\partial u_h}{\partial x_k}(x)\frac{\partial^2 g_{ij}}
{\partial x_h \partial x_l}(x) +
\frac{\partial u_h}{\partial x_l}(x)\frac{\partial^2 g_{ij}}
{\partial x_h \partial x_k}(x)\, + \right.$$
$$\hspace*{2.3cm}
\left. \frac{\partial u_h}{\partial x_i}(x)\frac{\partial^2 g_{hj}}
{\partial x_k \partial x_l}(x)+
\frac{\partial u_h}{\partial x_j}(x)\frac{\partial^2 g_{hi}}
{\partial x_k \partial x_l}(x)\right)\,.$$
Permuting the indices $i,k$, in (16), and subtracting we obtain
$$g_{ii}(x)\frac{\partial^3 u_i}{\partial x_j \partial x_k \partial
x_l}(x) -
g_{kk}(x)\frac{\partial^3 u_k}{\partial x_j \partial x_i \partial
x_l}(x) = \lambda_{kjil} - \lambda_{ijkl}\,,$$
and permuting the indices $j,k$,
\begin{equation}
g_{ii}(x)\frac{\partial^3 u_i}{\partial x_k \partial x_j \partial
x_l}(x) -
g_{jj}(x)\frac{\partial^3 u_j}{\partial x_k \partial x_i \partial
x_l}(x) = \lambda_{jkil} - \lambda_{ikjl}\,.
\end{equation}
By adding (16) and (17), we obtain
$$2g_{ii}(x)\frac{\partial^3 u_i}{\partial x_j \partial x_k \partial x_l}(x)=
\lambda_{jkil} - \lambda_{ikjl} - \lambda_{ijkl}\,.$$
Formula (7) then yields\newpage
\begin{equation}
g_{ii}(x)\frac{\partial^3 u_i}{\partial x_j \partial x_k \partial x_l}(x)=
\frac{1}{3}\sum_{h=1}^n\left[
\frac{\partial u_h}{\partial x_i}(x)\left(R_{klhj}(x)+R_{khlj}(x)\right)+
\right.
\end{equation}
$$\hspace*{5cm}\frac{\partial u_h}{\partial x_j}(x)\left(R_{ikhl}(x)
+R_{ihkl}(x)\right)+$$
$$\hspace*{5cm}\frac{\partial u_h}{\partial x_k}(x)
\left(R_{ihjl}(x)+R_{ijhl}(x)\right)+$$
$$\hspace*{5cm}\left.
\frac{\partial u_h}{\partial x_l}(x)
\left(R_{kihj}(x)+R_{khij}(x)\right)\right].$$
Since the left hand side of (18) is symmetric with respect to the indices
$j, k, l$, permuting $j$ and $l$, and equating the corresponding right
hand sides we obtain (15), thus also showing that the right hand side of
(18) vanishes and hence (18) reduces to (14). The proof is thus complete.
\fin

\vspace{2mm}

\noindent {\em Remark 14.} Equation (15) is invariant under the group of
order 8 generated by the permutations $\gamma:(ijkl)
\mapsto (jkli)$, and $\tau:(ijkl)\mapsto (ilkj)$. Note that $\gamma\com
\tau = \tau\com \gamma^2$. Accordingly, in examining (15) we only need
to consider the following three cases: $1^{st})\,\, i\leq j\leq k\leq l$,
$\;\;2^{nd})\,\, i\leq k< j\leq l$, $\;\;3^{rd})\,\, i\leq j\leq l< k$.

\vspace{2mm}

\noindent {\em Theorem 15.} {\bf (i)} If dim$N = n =1$, for every
$j^2_x g\in J^2({\cal M})$, $\Phi^2_{j^2_x g}$ is bijective.

\vspace{1mm}

\noindent {\bf (ii)} If dim$N = n = 2$, for every $j^2_x g\in
J^2({\cal M})$, rk.$\Phi^2_{j^2_x g} = 19$.

\vspace{1mm}

\noindent {\bf (iii)} For each $n\geq 3$, there exists a dense open subset
${\cal O}^{n,2}\subset J^2({\cal M})$, such that
\hspace*{7mm} for every $j^2_x g\in {\cal O}^{n,2}$, $\Phi^2_{j^2_x g}$
is injective.

\vspace{2mm}

\noindent {\em Proof.} {\bf (i)} From (8) and (14) it follows that
$\Phi^2$ is injective. Hence rk.$\Phi^2_{j^2_x g} =$ dim$J^3_x(TN) = 4 =$
dim$T_{j^2_x g}(J^2 N)$.\newline
{\bf (ii)} Using the above remark it is not difficult
to check that equation (15) is identically satisfied if $n=2$. Accordingly,
from (8) and (14) it follows that a $3$--jet $j^3_x X\in
\hbox{\rm Ker}\Phi^2_{j^2_x g}$ is completely determined by
$(\partial u_2/\partial x_1)(x)$. Hence the kernel
of $\Phi^2$ is a vector subbundle of rank $1$, and therefore
rk.$\Phi^2_{j^2_x g} =$ dim$J^3_x(TN) - 1 = 20 - 1 = 19$.\newline
{\bf (iii)} Let $r = \sum_{i,j} r_{ij} dx_i\otimes dx_j$ be the Ricci tensor
of $g$; {\em i.e.}, $r(X,Y) = \hbox{\rm trace of}\,Z\mapsto R(Z,X)(Y)$.
Then we have: $r_{ij}(x) = \sum_h g_{hh}(x) R_{hjhi}(x)$. Since $r$ is
symmetric there exists a unique endomorphism $A:T_x(N)\rightarrow T_xN$,
such that for every $X,Y\in T_x(N)$,
\begin{equation}
r(X,Y) = g(A(X),Y) = g(X,A(Y))\,.
\end{equation}
Moreover,
let $B$ be the endomorphism given by: $B = \sum_{i,j}(\partial u_i/
\partial x_j)(x) d_x x_j\otimes (\partial/\partial x_i)_x$. From the
second equation of (8) we deduce that for every $X,Y\in T_x(N)$,
\begin{equation}
g(B(X),Y) + g(X,B(Y)) = 0\,.
\end{equation}
Let ${\cal O}^{n,2}$ be the set of points
$j^2_x g$ such that the eigenvalues of $A$ in $T_x(N)\otimes
\hbox{\Bbb C}$ are pairwise different. We shall prove that on
${\cal O}^{n,2}$ the unique solution of (15) is the trivial one.
Let us denote by $<X,Y>$ the bilinear form induced by $g_x$ on the
complex vector space $T_x(N)\otimes \hbox{\Bbb C}$, so that (19) and
(20) imply
$$<A(Z),W> = <Z, A(W)>\hspace*{1cm} <B(Z),W> + <Z, B(W)> = 0\,,$$
for every $Z,W\in T_x(N)\otimes \hbox{\Bbb C}$.
Letting $h=k$ in (15), multiplying by $g_{kk}(x)$, and using the second
equation of (8), we obtain
$$\sum_{h=1}^n \left(r_{jh}(x)\frac{\partial u_h}{\partial x_i}(x) +
r_{ih}(x)\frac{\partial u_h}{\partial x_j}(x)\right) = 0\,,$$
or equivalently,
$$<AB(Z), W> + <Z, AB(W)> = 0\,,$$
for every $Z,W\in T_x(N)\otimes \hbox{\Bbb C}$.
Let $A(Z_i) = \lambda_i\,Z_i$ be the eigenvalues (and the eigenvectors) of
$A$. From the above equations we then obtain: $<AB(Z_i), Z_j>+<Z_i,AB(Z_j)>
=0$, or equivalently, $<B(Z_i),A(Z_j)>+<A(Z_i),B(Z_j)>=0$; {\em i.e.},
$\lambda_j <B(Z_i),Z_j>+ \lambda_i<Z_i, B(Z_j)>=0$. Hence $(\lambda_i
-\lambda_j)<Z_i,B(Z_j)>=0$. By virtue of the hypothesis this implies
$B(Z_j)=0$, and since $Z_1,...,Z_n$ is a basis of the complex tangent
space, we have $B=0$.
\fin

\vspace{8mm}

\noindent {\bf 7. The rank of $\Phi^r$, $r\geq 3$}

\vspace{5mm}

\noindent {\em Theorem 16.} For each $r\geq 3$, there exists a dense open
subset ${\cal O}^{n,r}\subset J^r({\cal M})$, such that for every $j^r_x g\in
{\cal O}^{n,r}$, $\Phi^r_{j^r_x g}$ is injective.

\vspace{2mm}

\noindent {\em Proof.} First we prove the theorem for $r=3$. We distinguish
two cases:

\vspace{1mm}

\noindent {\bf (i)} dim$N = n \neq 2$. For the sake of simplicity we write
${\cal O}^{1,2} = J^2({\cal M})$. From (i) and (iii) in Theorem 15 we know that
$\Phi^2$ is injective on ${\cal O}^{n,2}$. We set ${\cal O}^{n,3} = p_{32}
^{-1}({\cal O}^{n,2})$. Assume $j^3_x g\in {\cal O}^{n,3}$. Since
$\Phi^2_{|{\cal O}^{n,2}}$ is injective from formula (3) we have that:
$j^4_x(X)\in$Ker$\Phi^3_{j^3_x g}$ if and only if $j^3_x(X) = 0$, and
furthermore for every $i,j,k,l,m=1,...,n$,
\begin{equation}
\frac{\partial^4 u_j}{\partial x_i \partial x_k \partial x_l \partial
x_m}(x) g_{jj}(x) +
\frac{\partial^4 u_i}{\partial x_j \partial x_k \partial x_l \partial
x_m}(x) g_{ii}(x) = 0\,.
\end{equation}
Permuting $i$ and $k$ in (21) and subtracting we have
$$g_{ii}(x) \frac{\partial^4 u_i}{\partial x_j \partial x_k \partial x_l
\partial x_m}(x) - g_{ii}(x) \frac{\partial^4 u_i}{\partial x_j \partial
x_k \partial x_l \partial x_m}(x) = 0\,,$$
and permuting the indices $j,k$, and adding the equation thus obtained to
(21) we have
$$\frac{\partial^4 u_i}{\partial x_j \partial x_k \partial x_l
\partial x_m}(x) = 0\,,\,\hbox{\rm for every}\;\, i,j,k,l,m=1,...,n\,.$$
Hence $j^4_x(X) = 0$, and accordingly $\Phi^3_{|{\cal O}^{n,3}}$ is injective.

\vspace{1mm}

\noindent{\bf (ii)} dim$N = n = 2$. From formula (3) we conclude that
$j^4_x(X)$ belongs to the kernel of $\Phi^3_{j^3_x g}$, if and only if
in addition to equations (8) and (14) the following conditions hold true:
For every $|\alpha| = 3\,,\,i,j=1,2$,
$$ \sum_{h=1}^2\left(\frac{\partial u_h}{\partial x_i}(x)\frac{\partial^3
g_{hj}}{\partial x^\alpha}(x) +
\frac{\partial u_h}{\partial x_j}(x)\frac{\partial^3
g_{hi}}{\partial x^\alpha}(x) +
\sum_{k=1}^2 \alpha_k \frac{\partial u_h}{\partial x_k}(x)
\frac{\partial^3 g_{ij}}{\partial x^{\alpha-(k)+(h)}}(x)\right) + $$
\begin{equation}
g_{jj}(x) \frac{\partial^4 u_j}{\partial x^{\alpha+(i)}}(x) +
g_{ii}(x) \frac{\partial^4 u_i}{\partial x^{\alpha+(j)}}(x) = 0\,.
\end{equation}
Recall that equation (15) is identically satisfied if $n = 2$. Moreover
the expansion given in (5) yields (cf. [3]):
\begin{eqnarray*}
\frac{\partial^3 g_{ij}}{\partial x_k \partial x_l \partial x_m}(x)=
\frac{1}{6}\frac{\partial}{\partial x_k}(R_{imjl}+R_{iljm})(x) & + &
\frac{1}{6} \frac{\partial}{\partial x_l}(R_{imjk}+R_{ikjm})(x)  \\
& + & \frac{1}{6} \frac{\partial}{\partial x_m}
(R_{iljk}+R_{ikjl})(x)\,.
\end{eqnarray*}
{}From the above equation and (22) we then obtain:
\begin{eqnarray*}
& & (\alpha=(3,0),i=j=1):\;\;\frac{\partial^4 u_1}{\partial x_1^4}(x)=0, \\
& & (\alpha=(3,0),i=1,j=2):\;\;g_{22}(x)\frac{\partial^4 u_2}{\partial
x_1^4}(x) + g_{11}(x)\frac{\partial^4 u_1}{\partial x_1^3 \partial
x_2}(x)=0, \\
& & (\alpha=(3,0),i=j=2):\;\;\frac{\partial u_2}{\partial x_1}(x)
\frac{\partial R_{1212}}{\partial x_2}(x) + 2g_{22}(x) \frac{\partial^4
u_2}{\partial x_1^3 \partial x_2}(x)=0, \\
& & (\alpha=(2,1),i=j=1):\;\;\frac{\partial^4 u_1}{\partial x_1^3
\partial x_2}(x)=0,  \\
& & (\alpha=(2,1),i=1,j=2):\;\;-\frac{1}{3}\frac{\partial u_2}{\partial
x_1}(x)\frac{\partial R_{1212}}{\partial x_2}(x) + g_{11}(x)\frac{\partial
^4 u_1}{\partial x_1^2 \partial x_2^2}(x) + \\
& & \hspace*{8.3cm} g_{22}(x)\frac{\partial^4 u_2}{\partial x_1^3
\partial x_2}(x)=0, \\
& & (\alpha=(2,1),i=j=2):\;\;-\frac{1}{3}\frac{\partial u_2}{\partial x_1}
(x)\frac{\partial R_{1212}}{\partial x_1}(x)+2g_{11}(x)\frac{\partial^4
u_2}{\partial x_1^2 \partial x_2^2}(x)=0, \\
& & (\alpha=(1,2),i=j=1):\;\;\frac{1}{3}\frac{\partial u_2}{\partial x_1}
(x)\frac{\partial R_{1212}}{\partial x_2}(x)+2g_{11}(x)\frac{\partial^4
u_1}{\partial x_1^2 \partial x_2^2}(x)=0, \\
& & (\alpha=(1,2),i=1,j=2):\;\;\frac{1}{3}g_{22}(x)
\frac{\partial u_2}{\partial x_1}(x)\frac{\partial R_{1212}}{\partial x_1}(x)
+ \frac{\partial^4 u_1}{\partial x_1 \partial x_2^3}(x) + \\
& & \hspace*{7.3cm} g_{11}(x)g_{22}(x)\frac{\partial^4 u_2}{\partial x_1^2
\partial x_2^2}(x)=0, \\
& & (\alpha=(1,2),i=j=2):\;\;\frac{\partial^4 u_2}{\partial x_1
\partial x_2^3}(x)=0, \\
& & (\alpha=(0,3),i=j=1):\;\;-g_{22}(x)\frac{\partial u_2}{\partial x_1}
\frac{\partial R_{1212}}{\partial x_1}(x)+2\frac{\partial^4
u_1}{\partial x_1 \partial x_2^3}(x)=0, \\
& & (\alpha=(0,3),i=1,j=2):\;\;g_{22}(x)\frac{\partial^4 u_2}{\partial x_1
\partial x_2^3}(x)+g_{11}(x)\frac{\partial^4 u_1}{\partial x_2^4}(x)=0,\\
& & (\alpha=(0,3),i=j=2):\;\;\frac{\partial^4 u_2}{\partial x_2^4}(x)=0.
\end{eqnarray*}
It is not difficult to check that the above system is equivalent to saying
that for every $i,j,k,l,m=1,2$,
$$\frac{\partial^4 u_i}{\partial x_j \partial x_k \partial x_l \partial
x_l}(x) = 0\,,$$
and
$$\frac{\partial u_2}{\partial x_1}(x)\frac{\partial R_{1212}}{\partial
x_1}(x) = 0\hspace*{1cm} \frac{\partial u_2}{\partial x_1}(x)
\frac{\partial R_{1212}}{\partial x_2}(x) = 0\;.$$
Hence in the dense open subset ${\cal O}^{2,3}$ of the $2$--jets of metrics
of a surface whose curvature satisfies  $\parallel (\nabla R)(x)\parallel > 0$,
we have $j^4_x(X) = 0$, and therefore $\Phi^3_{|{\cal O}^{2,3}}$ is injective.

By induction on $r$ we now prove the general statement of the theorem. For
every $r\geq 3$, $n\geq 1$, we set ${\cal O}^{n,r}=p^{-1}_{r 3}({\cal O}
^{n,3})$. Assume $j^{r+1}_x(X)\in Ker\Phi^r_{j^r_x g}$, with $j^r_x g\in
{\cal O}^{n,r}$. Since $\overline{X}^r_{j^r_x g}$ projects onto
$\overline{X}^{r-1}_{j^{r-1}_x g}$, it follows from the induction hypothesis
that $j^r_x(X) = 0$, and from formula (3) we thus deduce
$$\frac{\partial^{r+1}u_j}{\partial x^{\alpha+(i)}}(x)g_{jj}(x)+
\frac{\partial^{r+1}u_i}{\partial x^{\alpha+(j)}}(x)g_{ii}(x) = 0\,,$$
for every $i,j=1,...,n$, $|\alpha|=r$. Let $k$ be an index
such that $\alpha_k > 0$. We set $\beta =\alpha - (k)$,
so that the above equation reads as follows
\begin{equation}
g_{ii}(x)\frac{\partial^{r+1}u_i}{\partial x_j\partial x_k \partial
x^\beta}(x)+
g_{jj}(x)\frac{\partial^{r+1}u_j}{\partial x_i\partial x_k\partial x
^\beta}(x)= 0\,.
\end{equation}
Permuting the indices $i,k$, in (23), and subtracting we have
$$g_{ii}(x)\frac{\partial^{r+1}u_i}{\partial x_j\partial x_k \partial
x^\beta}(x) -
g_{kk}(x)\frac{\partial^{r+1}u_k}{\partial x_i\partial x_j\partial x
^\beta}(x)= 0\,,$$
and permuting the indices $j,k$, and adding the above equation to (23),
we obtain
$$\frac{\partial^{r+1}u_i}{\partial x_j\partial x_k \partial
x^\beta}(x)=0\,,$$
thus proving that $j^{r+1}_x(X) = 0$, and finishing the proof of the theorem.
\fin

\vspace{8mm}

\noindent {\bf 8. Calculating $i_{n,r}$}

\vspace{5mm}

\noindent {\em Theorem 17.} On a dense open subset of $J^r({\cal M})$, the
number $i_{n,r}$ of functionally independent metric differential invariants
is the following:
\begin{description}
\item[(i)] for every $n\geq 1$, $i_{n,0} = i_{n,1} = 0$,
\item[(ii)] for every $r\geq 0$, $i_{1,r} = 0$,
\item[(iii)] $i_{2,2} = 1$, and for every $r\geq 3$, $i_{2,r} = \frac{1}
{2} (r+1)(r-2)$,
\item[(iv)] for every $n\geq 3$, $r\geq 2$,
$$i_{n,r} = n + \frac{(r-1)n^2-(r+1)n}{2(r+1)}{n+r\choose r}\,.$$
\end{description}

\vspace{2mm}

\noindent{\em Proof.} First of all we confine ourselves to the dense
open subset ${\cal O}^r$ prescribed in Theorem 10--(iii) where we know that
the metric differential invariants of order $r$ coincide with the ring of
first integrals of the fundamental distribution.

\vspace{1mm}

\noindent {\bf (i)} It follows from formula (2) that the
vector fields $\overline{X}_{j^0_x g}$ span $T_{j^0_x g}({\cal M})$; hence
$i_{n,0} = 0$. Since $\Phi^1$ is surjective (Theorem 12) from Corollary
11 we conclude that $i_{n,1} = 0$.

\vspace{1mm}

\noindent {\bf (ii)} From Theorem 16 we know that $\Phi^r$ is injective
(in fact ${\cal O}^{1,r} = J^r({\cal M})$ in this case). From Corollary 11 we
thus have $i_{1,r} =$ dim$J^r({\cal M}) -$ rk.$\Phi^r_{j^r_x g} =$
dim$J^r({\cal M}) -$ dim$J^{r+1}_x(TN) = (r+2) - (r+2) = 0$.

\vspace{1mm}

\noindent {\bf (iii)} That $i_{2,2} = 1$, follows directly from part
(ii) of Theorem 15. Moreover, the formula for $r\geq 3$, is a particular
case of the formula in (iv).

\vspace{1mm}

\noindent {\bf (iv)} From Theorem 15--(iii), Theorem 16 and Corollary
11 we have
\begin{eqnarray*}
i_{n,r} & = & \hbox{\rm dim}J^r({\cal M}) - \hbox{\rm dim}J^{r+1}_x(TN)  \\
& = & \left( n + \frac{n(n+1)}{2}{n+r\choose r}\right) - n{n+r+1\choose
r+1} \\
& = &  n + \frac{(r-1)n^2-(r+1)n}{2(r+1)}{n+r\choose r}\,.
\end{eqnarray*}
\fin

\vspace{2mm}

\noindent {\em Remark 18.} For $n\geq 3$, there is a classical procedure
in order to obtain second order metric invariants, the so--called
curvature invariants ([9], p. 146). In the generic case there is an
essentially unique frame reducing $g$ and its Ricci tensor to a canonical
form. The invariants are the components of the Weyl tensor on that frame
plus the $n$ eigenvalues of the Ricci tensor. Let us calculate the
dimension of the space of Weyl tensors. Following the same notations as
[2], 1.105--116, we have that the space ${\cal C}E$ of curvature tensors
(here $E = T^*_x(N)$), breaks into three irreducible subspaces under
the natural action of the orthogonal group, ${\cal C}E = {\cal U}E
\oplus {\cal Z}E \oplus {\cal W}E$, where dim ${\cal U}E = 1$, dim
${\cal Z}E = \frac{n(n+1)}{2} -1$ (traceless symmetric $2$--tensors),
and ${\cal W}E$ are the Weyl tensors. Hence:\newpage
\begin{eqnarray*}
\hbox{\rm dim}\, {\cal W}E & = & \hbox{\rm dim}\, {\cal C}E -
\frac{n(n+1)}{2} \\
& = & \hbox{\rm dim}\,S^2\bigwedge^2 E - \hbox{\rm dim}\, \bigwedge^4 E -
\frac{n(n+1)}{2} \\
& = & \frac{n^4 - 7n^2 - 6n}{2}\,.
\end{eqnarray*}
Hence the number of curvature invariants is
$$n + \frac{n^4 - 7n^2 - 6n}{2}  = \frac{n^4 - 7n^2 + 6n}{2} = i_{n,2}\,.$$
Accordingly, this shows that the number of functionally independent
curvature invariants are exactly the number of functionally independent
second order metric invariants.

\end{document}